\documentclass{SCIYA2017enOL}

\DeclareMathOperator*{\res}{Res}
\theoremstyle{plain}

\newtheorem{Assumption}{Assumption}
\theoremstyle{remark}

\online
\begin{document}

\ensubject{fdsfd}

\ArticleType{ARTICLES}
\Year{2023}
\Month{January}%
\Vol{66}
\No{1}
\BeginPage{1} %
\DOI{10.1007/s11425-016-5135-4}
\ReceiveDate{January 1, 2022}
\AcceptDate{January 1, 2022}
\OnlineDate{January 1, 2022}

\title[]{Long time asymptotic  behavior for  the nonlocal nonlinear  Schr\"odinger equation in the  solitonic region}
{Long time asymptotic  behavior for  the nonlocal nonlinear  Schr\"odinger equation in the  solitonic region}

\author[]{Gaozhan Li}{gzli20@fudan.edu.cn}
\author[]{Yiling YANG}{ylyang19@fudan.edu.cn}
\author[$\ast$]{Engui Fan}{faneg@fudan.edu.cn}

\AuthorMark{Li G. Z.}

\AuthorCitation{Li G. Z., Yang Y. L., Fan E. G.}

\address[1]{School of Mathematical Sciences  and Key Laboratory   for Nonlinear Science, Fudan   University, Shanghai 200433, P. R. China.}

\abstract{In this paper, we  extend  $\overline\partial$-steepest descent method  to study the Cauchy problem for the nonlocal nonlinear
	Schr\"odinger (NNLS)  equation with  weighted Sobolev  initial data
	\begin{align*}
		&{\rm i}q_{t}+q_{xx}+2\sigma q^2(x,t)\overline{q}(-x,t)=0, \\
		& q(x,0)=q_0(x),
	\end{align*}	
	where  $ q_0 \in H^{1,1}(\mathbb{R}) $.
	Based on  the spectral analysis of the Lax pair,  the  solution of the   Cauchy problem is expressed
	in terms of  the  solution   of  a Riemann-Hilbert problem,  which is transformed into a solvable model after a series of deformations.
	We further obtain the  asymptotic expansion of  the solution  to  the Cauchy problem for  the  NNLS equation in the   solitonic region.
	The leading   term is soliton solutions, the second term is     the interaction between solitons and dispersion,
	the error term  comes from a corresponding $\bar{\partial}$-problem.
	Comparing to the asymptotic results on the  classical NLS equation,
	the major difference is  the second and third terms in asymptotic expansion for the
	NNLS equation  were affected by a  function  depends on the scattering data and  the stationary phase point.}



\maketitle

\section{Introduction}

\quad
In this paper, we  extend  $\overline\partial$-steepest descent method  to study the Cauchy problem for the nonlocal
Schr\"odinger (NNLS)  equation with  weighted Sobolev  initial data
\begin{align}
	&{\rm i}q_{t}+q_{xx}+2\sigma q^2(x,t)\overline{q}(-x,t)=0, \label{NNLS} \\
	& q(x,0)=q_0(x),\label{cauchy}
\end{align}	
where  $ \sigma=\pm1 $  and $ q_0\in H^{1,1}(\mathbb{R})$, which is a  weighted Sobolev space  defined by
$$H^{1,1}(\mathbb{R})=\left\{f\in L^2(\mathbb{R}) :   f_x\in L^2(\mathbb{R}),    (1+\left<\cdot\right> ) f \in L^2 (\mathbb{R})\right\},
$$
where $\langle x\rangle=\sqrt{1+x^2}$.
The  NNLS  equation  was first proposed
by Ablowitz and Musslimani \cite{ablo1,ablo2} and has attracted much
attention in recent years due to its distinctive properties.
It was   shown that the  NNLS equation   is invariant under the joint transformations $x  \to -x$, $t \to -t$ and complex conjugation \cite{PT1,PT2}.
PT symmetric and non-Hermitian physics has been
the subject of an intense research for the last decade most notably
in classical optics, quantum mechanics and topological photonics \cite{13,14,15,111,112,113,pt1,add1}.  PT symmetry of the NNLS equation amounts to the invariance of
the so-called self-induced potential in the case of classical optics \cite{PT3}.   Wave propagation in PT
symmetric coupled waveguide or photonic lattices also has been observed experimentally.
In recent years, much work    on  the mathematical structure and properties for the
NNLS equation   also have been widely studied   \cite{ablo3,ablo4,28}.

The inverse scattering transform (IST) method was first used to solve the NNLS equation by Ablowitz and  Musslimani  \cite{ablo1}.
The long-time asymptotics  for the  NNLS equation with rapidly decaying initial data  was first considered  due to Rybalko and Shepelsky
with  nonlinear steepest decent method \cite{jmp-z}.  In their series of papers, they further
studied long-time asymptotics  for the  NNLS equation with step-like initial data \cite{steplike1,steplike2}.
The major  difference from classical NLS equation is that   the NNLS equation has different symmetries  on their Jost functions and
scattering data. Also the NNLS equation admits two  scattering  coefficients   $r_1(k)$ and $r_2(k)$, which result in   $1+\sigma r_1(k)r_2(k)$ is, in general,  not real-valued function.
Therefore the   scattering data   ${\rm Im} \nu(-\xi)$ (\ref{nu})  also   affect  long time asymptotic results.
These  bring about    difficulties  to  the classify   discrete spectrum and analyze asymptotics of solutions.
The nonlinear steepest descent method    was first  developed  by  Deift and  Zhou to   investigate the long-time behavior of mKdV equation \cite{DZAnn}.
Later this method is further   generalized into  a powerful   $\bar{\partial}$-steepest descent method  to analyze asymptotic of orthogonal polynomials
by McLaughlin and Miller \cite{M&M2006, M&M2008}. In recent years, this  $\bar{\partial}$-method
has been widely used to  analyze  long-time asymptotics of  integrable systems \cite{Dieng,dbar1,Liu3,YF1}.

In this paper,  we would like to extend above $\bar{\partial}$-method to
analyze  long time asymptotics of the NNLS equation with weighted Sobolev initial data.
Compared with work  \cite{jmp-z},  we consider more general   weighted Sobolev initial data which support
appearance of solitions.  And because of the more general initial data the  ${\rm Im}\nu(-\xi) $ comes from reflect coefficient  gives a less precise estimate  to the asymptotic expansion of the NNLS equation.
Moreover, to our knowledge  $\bar\partial$-steepest decent method  has not been applied to any nonlocal integrable system.

This paper is arranged as follows.  In Section \ref{sec2},
we recall  some main  results  in  the construction  process  of the  RH  problem with respect to  the initial problem   (\ref{NNLS})-(\ref{cauchy})
obtained in  \cite{jmp-z}.   We further prove that the scattering coefficients $r_1(k)$ and $ r_2(k)$ belong to the  Sobolev space $H^1(\mathbb{R})$.
In Section \ref{sec3},  we   transform the RH problem   into a mixed $\overline{\partial}$-RH problem, which is  solved    by separating
it into a pure $\overline{\partial}$-problem and a pure RH problem.
The  pure RH problem  can be estimated by a solvable model RH problem  and  a soliton model   RH  problem.
The pure $\overline{\partial}$-problem is  analyzed  for large  $t\to\infty$.
Finally, we obtain  long time asymptotic  behavior of   the NNLS equation   by using the reconstruction formula  in Section \ref{sec4}.

\section {Direct and inverse scattering transform }\label{sec2}

\subsection{The Jost functions }

\quad The NNLS equation (\ref{NNLS})    admits the Lax pair \cite{jmp-z}
\begin{equation}
	\Phi_x = (-{\rm i}k\sigma_3+Q) \Phi,\hspace{0.5cm}\Phi_t =(-2{\rm i}k^2\sigma_3+V)\Phi, \label{lax0}
\end{equation}
where
\begin{align}
	&Q=\left(\begin{array}{cc}
		0 & q(x,t) \\
		-\sigma \overline{q(-x,t)} & 0
	\end{array}\right),
	V=\left(\begin{array}{cc}
		A & B \\
		C & -A
	\end{array}\right),\nonumber\\
	&A={\rm i}\sigma q(x,t)\overline{q(-x,t)}, \ \ \
	B=2kq(x,t)+{\rm i}q_{x}(x,t), \nonumber\\
	&C=-2k\sigma \overline{q(-x,t)}+{\rm i}\sigma (\overline{q(-x,t)})_x.\nonumber
\end{align}
Here  $k\in \mathbb{C}$ is spectral parameter, and   $ \sigma_3$ is the third Pauli matrix.

Under the  initial value  (\ref{cauchy}),    the Lax pair (\ref{lax0})  have matrix-valued Jost solutions
$$\Phi_j=\Psi_j {\rm e}^{(-{\rm i}kx-2{\rm i}k^2 t)\sigma_3}, \ j=1, 2, $$
where  $ \Psi_j$ solve following Volterra integral equation
\begin{align}
	&\Psi_1(x,t;k)=I+\int_{-\infty}^{x}{\rm e}^{{\rm i}k(y-x)\sigma_3}Q(y,t)\Psi_1(y,t;k){\rm e}^{-ik(y-x)\sigma_3}dy,\label{volltra1}\\	
	&\Psi_2(x,t;k)=I+\int_{+\infty}^{x}{\rm e}^{{\rm i}k(y-x)\sigma_3}Q(y,t)\Psi_2(y,t;k){\rm e}^{-{\rm i}k(y-x)\sigma_3}dy.\label{volltra2}
\end{align}
It was shown  that  for $ q(\cdot,t) \in L^{1}(\mathbb{R})$ for $t\geq0$,     the first column of $ \Psi_1 $ and second  columns of $ \Psi_2 $ are analytic in
$ k\in \mathbb{C}^+$ and the first column of $ \Psi_2 $ and second  columns of $ \Psi_1 $ are analytic in $ k\in \mathbb{C}^-$, respectively \cite{jmp-z}.   Here, $ \mathbb{C}^{\pm}=\left\lbrace  z\in\mathbb{C}|\pm{\rm Im}(z)>0 \right\rbrace $.
There exists  a scattering matrix $ S(k)$  such that
\begin{equation}
	\Phi_1(x,t;k)=\Phi_2(x,t;k)S(k),\hspace{0.5cm}S(k)=\left(\begin{array}{cc}
		a_1(k)&\tilde{b}(k)\\
		b(k)&a_2(k)\\
	\end{array}\right). \label{Smtx}
\end{equation}
$\Phi_1(x,t;k)$ and $ \Phi_2(x,t;k)$ admit the following symmetry
\begin{equation}
	\Lambda \overline{\Phi_1(-x,t,-k)} \Lambda^{-1}=\Phi_2(x,t;k),\  \Lambda=\left(\begin{array}{cc}
		0&\sigma \\
		1&0\\
	\end{array}\right).\label{sym}
\end{equation}
Above symmetry  condition can be expended as
\begin{align*}
	\begin{pmatrix}
		0&\sigma\\
		1&0\\
	\end{pmatrix}\overline{[\Phi_1]_1(-x,-\overline{k})}=
	[\Phi_2]_2(x,k), k\in\mathbb{C}^+,\\
	\begin{pmatrix}
		0&1\\
		\sigma&0\\
	\end{pmatrix}\overline{[\Phi_1]_2(-x,-\overline{k})}=
	[\Phi_2]_1(x,k), k\in\mathbb{C}^-,\\
\end{align*}where $[\cdot]_j$, $j=1,2$ referring the $j$th column of matrix.
Combining with the (\ref{Smtx}), we have
\begin{align}
	& \tilde{b}(k)=-\sigma \overline{b(-k)},k\in\mathbb{R},\\
	&a_1(k)=\overline{a_1(-\overline{k})},\ k\in{\mathbb{C}^+}\cup\mathbb{R},\label{syma1}\\
	&a_2(k)=\overline{a_2(-\overline{k})},\ k\in{\mathbb{C}^-}\cup\mathbb{R}.\label{syma2}
\end{align}
Thus, we get the relationship of $ \Psi_j,\ (j=1,2)$
\begin{align}
	\Psi_1(x,t;k)&=\Psi_2(x,t;k){\rm e}^{-{\rm i}t\theta \hat{\sigma}_3}S(k),
\end{align}
where $\theta(x,t;k) =4k\xi+2k^2$ with  $\xi=\frac{x}{4t}\in\mathbb{R} $ is the phase function, and $\hat{\sigma}_3$ is a mark means for arbitrary $2\times2$ matrix $A,B$
\begin{equation*}
	{\rm e}^{\hat{A}}B={\rm e}^{A}B{\rm e}^{-A}.
\end{equation*}
We define a matrix-valued function
\begin{align}
	M(x,t;k)=\left\{ \begin{array}{ll}
		\left( \frac{[\Psi_1]_1(x,t;k)}{a_1(k)}, [\Psi_2(x,t;k)]_2\right),   & k\in \mathbb{C}^+,\\[12pt]
		\left( [\Psi_2(x,t;k)]_1,\  \frac{[\Psi_1]_2(x,t;k)}{a_2(k)}\right) , &k\in\mathbb{C}^-.\\
	\end{array}\right.
\end{align}
The $[\Psi_i]_j$ , $i,j=1,2$ means the first or second column of $\Psi_i$, respectively.
Moreover,  $r_1(k)$ and $ r_2(k)$ are the reflection coefficients defined by
\begin{equation}
	r_1(k)=\frac{b(k)}{a_1(k)},\ \ r_2(k)=\frac{\overline{b(-k)}}{a_2(k)}.
\end{equation}
The zeros of  $a_1(z),a_2(z)$ on $\mathbb{R}$   are known to
occur and they correspond to spectral singularities.  They are excluded from our analysis in this paper.  To avoid trouble from varieties of  kinds of poles and make $r_1(k), r_2(k)$ have enough smoothness and decaying property,
we discuss under following assumption:
\begin{Assumption}\label{asmp1}
	$a_1(k)$ in  $\mathbb{C}^+$  and $ a_2(k) $ in  $\mathbb{C}^-$  have the same amount of  zeros    respectively, which are   simple zeros  and  nonzero on the real axis $\mathbb{R}$.
\end{Assumption}

In this subsection, we establish  map from the initial data $q_0(x)$ to scattering coefficient $r_1(k),\ r_2(k)$.  From the general discussion from \cite{zhou1989}, to ensure the solvability  of the Riemann-Hilbert problem in the section \ref{RHP} and do the following $\bar{\partial}$-analysis, some  regularity of $r_1,r_2$ are required.  According to \cite{ZL2}, the following result is known, we list it here without proving.
\begin{proposition}
	When the initial value $q_0 \in H^{1,1}(\mathbb{R})$, the scattering coefficient $r_1,r_2\in H^1(\mathbb{R})$.
\end{proposition}

\subsection{A Riemann-Hilbert problem}\label{RHP}

Suppose that $a_1(k)$  has  simple zeros $\omega_1,...,\omega_{N_1}$ on $\{k\in\mathbb{C}|\ {\rm Im}k>0,\ {\rm Re}k>0\}$ and $\eta_1,..,\eta_{M_1}\in i\mathbb{R}^+$, while $a_2(k)$ has  simple zeros $\gamma_1,...,\gamma_{N_2}$ on $\{k\in\mathbb{C}|\ {\rm Im}k<0,\ {\rm Re}k>0\}$ and $\tau_1,..,\tau_{M_2}\in i\mathbb{R}^-$. Assumption \ref{asmp1} keeps the same number of zeros $2N_1+M_1=2N_2+M_2$ between the two half-plane.

The symmetry in (\ref{syma1}) and (\ref{syma2})
imply that
\begin{align*}
	&a_1(\omega_n)=0 \Longleftrightarrow \overline{a_1(- \bar{\omega}_n)}=0, \hspace{0.5cm}n=1,...,N_1,  \\
	&  a_2(\gamma_n)=0 \Longleftrightarrow \overline{a_2(- \bar{\gamma}_n)}=0,  \hspace{0.5cm}n=1,...,N_2.
\end{align*}
Therefore, the discrete spectrum is
\begin{align}
	\mathcal{N}_1&=\left\lbrace  \eta_m, \omega_n,-\overline{\omega_n}\ |m=1,..M_1, n=1,...N_1\right\rbrace,\\
	\mathcal{N}_2&=\left\lbrace  \tau_m, \gamma_n,-\overline{\gamma_m}\ |m=1,..M_2, n=1,...,N_2\right\rbrace.
\end{align}
For brief, denote \begin{equation*}
	-\overline{\omega_n}=\omega_{-n},n=1,...,N_1,\hspace{0.9cm}-\overline{\gamma_m}=\gamma_{-m},m=1,...,N_2.
\end{equation*}
Direct computation gives the residue conditions of $ M(x,t;k)$ as follows
\begin{align}
	\begin{aligned}
		\res_{k=\omega_n }M(x,t;k)&=\lim_{k\rightarrow \omega_n}M(x,t;k)\left(\begin{array}{cc}
			0 & 0 \\
			c_n{\rm e}^{2{\rm i}t\theta} & 0
		\end{array}\right), \\
		\res_{k=\eta_m }M(x,t;k)&=\lim_{k\rightarrow \eta_m}M(x,t;k)\left(\begin{array}{cc}
			0 & 0 \\
			e_m{\rm e}^{2{\rm i}t\theta} & 0
		\end{array}\right),
	\end{aligned}\label{res+}\\
	\begin{aligned}
		\res_{k=\gamma_n }M(x,t;k)&=\lim_{k\rightarrow \gamma_n}M(x,t;k)\left(\begin{array}{cc}
			0 & d_n{\rm e}^{-2{\rm i}t\theta} \\
			0 & 0
		\end{array}\right),\\
		\res_{k=\tau_m } M(x,t;k)&=\lim_{k\rightarrow \tau_m}M(x,t;k)\left(\begin{array}{cc}
			0 & f_m{\rm e}^{-2{\rm i}t\theta} \\
			0 & 0
		\end{array}\right).
	\end{aligned}\label{res-}
\end{align}
The norming constant $c_n,\ d_n,\ e_m,\ f_m$ which are independent of $x,t$ with:
\begin{align*}
	&a_1'(\omega_{n})c_n[\Phi_2]_2(x,t;\omega_n)=[\Phi_1]_1(x,t;\omega_n),\hspace{1.2cm}n=\pm1,...,\pm N_1,\\
	&a_2'(\gamma_{n})d_n[\Phi_2]_1(x,t;\gamma_n)=[\Phi_1]_2(x,t;\gamma_n),
	\hspace{1.2cm}n=\pm1,...,\pm N_2,
\end{align*}
and \begin{align*}
	&a_1'(\eta_{m})e_m[\Phi_2]_2(x,t;\eta_m)=[\Phi_1]_1(x,t;\eta_m),\hspace{1.2cm}m=1,...,M_1\\
	&a_2'(\tau_{m})f_m[\Phi_2]_1(x,t;\tau_m)=[\Phi_1]_2(x,t;\tau_m),\hspace{1.2cm}m=1,...,M_2.
\end{align*}
Denote above constants as sets
	\begin{align*}
		\mathcal{P}_1=&\left\{c_{n},e_{m}|n=\pm1,..\pm N_1,m=1,..M_1\right\},\\
		\mathcal{P}_2=&\left\{d_{n},f_{m}|n=\pm1,..\pm N_2,m=1,..M_2\right\}.
\end{align*}
From the symmetry \eqref{sym},
\begin{align}\label{symc}
	\begin{aligned}
		c_n\overline{c_{-n}}=-\frac{\sigma}{a_1'(\omega_{n})^2},\\
		d_n\overline{d_{-n}}=-\frac{\sigma}{a_2'(\gamma_{n})^2},
	\end{aligned}
\end{align}
and
\begin{align}\label{syme}
	\begin{aligned}
		|e_m|^2=\frac{\sigma}{|a_1'(\eta_m)|^2},\\
		|f_m|^2=\frac{\sigma}{|a_2'(\tau_m)|^2}.
	\end{aligned}
\end{align}
The expression is well-defined because of $0=a_1(\omega_n)=\det([\Psi_1]_1(x,t;\omega_n),[\Psi_2]_2(x,t;\omega_n))$ for any $n=\pm1,..,\pm N_1$ as an example.
Due to  above discussion, $ M(x,t;k)$ solves the following matrix-valued Riemann-Hilbert problem:

\noindent\textbf{RHP0\label{rhp0}}.  Find a matrix-valued function $M(x,t;k)$ which satisfies

$\blacktriangleright$ $M(x,t;k)$ is analytic in $\mathbb{C}\setminus\left( \mathbb{R} \cup\mathcal{N}_1\cup\mathcal{N}_2\right)$.

$\blacktriangleright$   $M(x,t;k)$ has continuous boundary values $M(x,t;k^\pm)$ on $\mathbb{R}$ with
\begin{equation}
	M(x,t;k^+)=M(x,t;k^-)J(x,t;k),
\end{equation}
where $ k^{\pm}$ means   limit of $ k\in \mathbb{R}$ and
\begin{equation}
	J(x,t;k)=\left(\begin{array}{cc}
		1+\sigma r_1(k) r_2(k) & \sigma r_2(k){\rm e}^{-2{\rm i}t\theta(x,t;k)} \\
		r_1(k){\rm e}^{2{\rm i}t\theta(x,t;k)} & 1
	\end{array}\right).
\end{equation}

$\blacktriangleright$ Residue condition: The residue of $ M(x,t;k)$ satisfies (\ref{res+}), (\ref{res-}).

$\blacktriangleright$ Asymptotic condition:
\begin{equation}
	M(x,t;k)=I+\mathcal{O}(k^{-1}),\hspace{0.5cm} k\rightarrow \infty.
\end{equation}
The reconstruction formula of $q(x,t)$ is given by
\begin{equation}
	q(x,t)=2{\rm i}\lim_{k\rightarrow\infty}(k[M(x,t;k)]_{12}).\label{qfml}
\end{equation}
Observing the RHP0, we can reconstruct the solution $q(x,t)$ by a set of scattering data
\begin{equation}\label{scat}
		\Big\lbrace r_1(k),r_2(k),\mathcal{N}_1,\mathcal{N}_2,\mathcal{P}_1,\mathcal{P}_2\Big\rbrace .
\end{equation}

\section {The deformation of the RH problem  }\label{sec3}

\quad Define \begin{equation}
	D^{\pm}=\left\lbrace z|\ {\rm Re}(z+\xi){\rm Im}(z+\xi)\in\mathbb{R}^{\pm}\right\rbrace,
\end{equation}
then for $k\in D^\pm,\ {\rm Re}(i\theta(x,t;k))\in \mathbb{R}^\pm $ as Figure \ref{D}.
Similar with \cite{jmp-z}, we give the assumption  of the scattering data as:
\begin{Assumption}\label{asym2}
	Let
	\begin{equation}
		\nu(k)=-\frac{1}{2\pi}\log(1+\sigma r_1(k)r_2(k)),\label{nu}
	\end{equation}
	with assuming
	$$| {\rm Im}\nu(k)|<\frac{1}{2},\ k\in\mathbb{R}.$$
\end{Assumption}
Introduce
\begin{equation}
	\delta(k)=\exp\left({\rm i}\int_{-\infty}^{-\xi}\frac{\nu(s)}{s-k}ds\right), k\in \mathbb{C}\setminus (-\infty,-\xi],
\end{equation}
then we prove that
\begin{figure}
	\centering
	\begin{tikzpicture}
		\draw[ ][-latex](-3,0)--(3,0);
		\draw[ ][-latex](0,-2)--(0,2);
		\fill (0.2,0) circle (0pt) node[below]{$-\xi$};
		\fill (1.2,1.2) circle (0pt) node[below]{$D_+$};
		\fill (-1.2,-1.2) circle (0pt) node[below]{$D_+$};
		\fill (1.2,-1.2) circle (0pt) node[below]{$D_-$};
		\fill (-1.2,1.2) circle (0pt) node[below]{$D_-$};
	\end{tikzpicture}
	\caption{  Denote  signature distribution  of  $ {\rm Re}(i\theta(x,t;k))$ as $D_\pm$.\label{D} }
\end{figure}
\begin{proposition}\label{dlt}
	$\delta(k) $ is analytic in $\mathbb{C}\setminus (-\infty,-\xi]$ and satisfies jump condition on $(-\infty,-\xi]$
	\begin{equation}
		\delta(k^+)=\delta(k^-)(1+\sigma r_1(k) r_2(k)).
	\end{equation}
	Moreover,
	\begin{equation}
		\delta(k)=1+\mathcal{O}(k^{-1}),\hspace{1.2cm} k\to\infty.
	\end{equation}
\end{proposition}
Noting that
\begin{equation}
	\chi(k)=-{\rm i}\nu(-\xi)\log(\xi+k+1)+{\rm i}\int_{-\xi-1}^{-\xi}\frac{\nu(s)-\nu(-\xi)}{s-k}+{\rm i}\int_{-\infty}^{-\xi-1}\frac{\nu(s)}{s-k},
\end{equation}
then a simple calculation gives
\begin{equation}
	\delta(k)=(\xi+k)^{{\rm i}\nu(-\xi)}\exp(\chi(k)),
\end{equation}
which  admits the following estimation
\begin{proposition}
	Let $\xi+k=r{\rm e}^{{\rm i}\phi},$ $|\phi|\leq \frac{\pi}{4},$ $r>0$, then
	\begin{equation}
		|\delta(k)-\delta_0 (\xi+k)^{{\rm i}\nu(-\xi)}|\lesssim |k+\xi|^{\frac{1}{2}-{\rm Im}\nu(-\xi)},\label{1/2}
	\end{equation}
	where $ \delta_0=\exp(\chi(-\xi))$.
\end{proposition}
\begin{proof}
	Simple calculation gives that	
	\begin{align}
		|\delta(k)-\delta_0 (\xi+k)^{{\rm i}\nu(-\xi)}|&=|(\xi+k)^{{\rm i}\nu(-\xi)}\exp(\chi(k))-\delta_0 (\xi+k)^{{\rm i}\nu(-\xi)}|\nonumber\\
		&\leq
		r^{-{\rm Im}\nu(-\xi)}{\rm e}^{-\phi {\rm Re}\nu(-\xi)}|\exp(\chi(k))-\exp(\chi(-\xi))|\nonumber\\
		\ &\lesssim r^{-{\rm Im}\nu(-\xi)}|\chi(k)-\chi(-\xi)|\nonumber\\
		&\lesssim |\xi+k|^{\frac{1}{2}-{\rm Im}\nu(-\xi)}.\nonumber
	\end{align}
	The proof of last line is similar with  that in  \cite{dbar1}.
\end{proof}

\subsection{A mixed $\bar\partial$-RH problem}
\quad
The long-time asymptotic  of RHP0  is affected by the growth and decay of the exponential function ${\rm e}^{\pm2{\rm i}t\theta}$ appearing in both the jump relation and the residue conditions. Therefore, we
define
\begin{align}
	\tilde{M}(x,t;k)=M(x,t;k)\delta(k)^{-\sigma_3}.
\end{align}
Through the Proposition \ref{dlt}, $	\tilde{M}(x,t;k)$ admits the following  Riemann-Hilbert problem

\noindent\textbf{RHP1}.  Find a matrix-valued function $\tilde{M}(x,t;k)$ which satisfies

$\blacktriangleright$ $\tilde{M}(x,t;k)$ is analytic in $\mathbb{C}\ \setminus\left( \  \mathbb{R}\  \cup\  \mathcal{N}_1 \ \cup\  \mathcal{N}_2\right) $.

$\blacktriangleright$ Jump condition: $\tilde{M}(x,t;k)$ has continuous boundary values $\tilde{M}(x,t;k^\pm)$ on $\mathbb{R}$ with
\begin{align}
	\tilde{M}(x,t;k^+)&=\tilde{M}(x,t;k^-)\tilde{J}(x,t;k),\ k\in \mathbb{R},\end{align}
where\footnotesize
\begin{align}
	&\tilde{J}(x,t;k) =\left\{ \begin{array}{ll}
		\left(\begin{array}{cc}
			1& \ \\
			\frac{r_1(k)}{1+\sigma r_1(k)r_2(k)}\delta^{-2}(k^-){\rm e}^{2{\rm i}t\theta}&1\\
		\end{array}\right)\left(\begin{array}{cc}
			1& \frac{\sigma r_2(k)}{1+\sigma r_1(k)r_2(k)}\delta^{2}(k^+){\rm e}^{-2{\rm i}t\theta} \\
			\ &1\\
		\end{array}\right),  & k<-\xi,\\[12pt]
		\left(\begin{array}{cc}
			1&\sigma r_2(k)\delta^2(k){\rm e}^{-2{\rm i}t\theta(x,t;k)} \\
			\ &1\\
		\end{array}\right)
		\left(\begin{array}{cc}
			1&\\
			r_1(k)\delta^{-2}(k){\rm e}^{2{\rm i}t\theta(x,t;k)}  &1\\
		\end{array}\right) , &k>-\xi.\\
	\end{array}\right.
\end{align}
\normalsize
$\blacktriangleright$ Residue condition: The residue of $ \tilde{M}(x,t;k)$ becomes
\begin{align}
	\begin{aligned}
		\res_{k=\omega_{n} }\tilde{M}(x,t;k)&=\lim_{k\rightarrow \omega_{n} }\tilde{M}(x,t;k)\left(\begin{array}{cc}
			0 & 0 \\
			\delta^{-2}(\omega_{n})c_n{\rm e}^{2{\rm i}t\theta} & 0
		\end{array}\right),\\
		\res_{k=\eta_{m} }\tilde{M}(x,t;k)&=\lim_{k\rightarrow \eta_{m}}\tilde{M}(x,t;k)\left(\begin{array}{cc}
			0 & 0 \\
			\delta^{-2}(\eta_{m})e_m{\rm e}^{2{\rm i}t\theta} & 0
		\end{array}\right),
	\end{aligned}\label{res1+}\\
	\begin{aligned}
		\res_{k=\gamma_{n} }\tilde{M}(x,t;k)&=\lim_{k\rightarrow \gamma_{n}}\tilde{M}(x,t;k)\left(\begin{array}{cc}
			0 & \delta^2(\gamma_{n}) d_n{\rm e}^{-2{\rm i}t\theta} \\
			0 & 0
		\end{array}\right),\\
		\res_{k=\tau_{m} }\tilde{M}(x,t;k)&=\lim_{k\rightarrow \tau_{m}}\tilde{M}(x,t;k)\left(\begin{array}{cc}
			0 & \delta^2(\tau_{m}) f_m{\rm e}^{-2{\rm i}t\theta} \\
			0 & 0
		\end{array}\right).
	\end{aligned}\label{res1-}
\end{align}

$\blacktriangleright$ Asymptotic condition:
\begin{equation}
	\tilde{M}(x,t;k)=I+\mathcal{O}(k^{-1}),\  k\rightarrow \infty.
\end{equation}

\quad Denote a positive number
\begin{equation}
	\rho_0=\frac{1}{2}\min\left\lbrace |x-y|>0\big\arrowvert x,y\in (\{- \xi\}\cup \mathcal{N}_1\cup \mathcal{N}_2)\right\rbrace .
\end{equation}
Let $\mathcal{Y}(k),\ k\in \mathbb{C} $ is a smooth function supported in $\left\lbrace k\big||k+\xi|<\frac{3}{2}\rho_0\right\rbrace $ and  $\mathcal{Y}=1 $ in  $\left\lbrace k\big||k+\xi|<\rho_0\right\rbrace $.
And $ \mathcal{X}(k),\ k\in \mathbb{C}$ is also  a smooth function supported in
$$\left\lbrace k\big||k-z|<\frac{3}{2}\rho_0,z\in \mathcal{N}_1\cup\mathcal{N}_2\right\rbrace $$
and $\mathcal{X}=1, \ k\in \left\lbrace k\big||k-z|<\rho_0,z\in \mathcal{N}_1\cup\mathcal{N}_2\right\rbrace $.
We hope to define the mixed $ \overline{\partial}$-RH problem, whose jump condition and $ \overline{\partial}$ derivative can be well controlled. Let
$$\Sigma_j =\left\lbrace k\Big| \xi+k=r{\rm e}^{{\rm i}\phi},\phi=\frac{2j-1}{4}\pi,r>0\right\rbrace,\ j=1,2,3,4,$$
which divides the complex plane $\mathbb{C}$ into six regions $\Omega_j, \ j=1, \cdots,6$  which are shown in Figure \ref{figR2}.
\begin{figure}[htp]
	\centering
	\begin{tikzpicture}[scale=0.7]
		
		\draw(0,0)--(3,3)node[right]{$\Sigma_1$};
		\draw(0,0)--(-3,3)node[left]{$\Sigma_2$};
		\draw(0,0)--(-3,-3)node[left]{$\Sigma_3$};
		\draw(0,0)--(3,-3)node[right]{$\Sigma_4$};
		\draw[dashed][->](-4,0)--(4,0);
		
		\draw[-latex](-3,-3)--(-1.5,-1.5);
		\draw[-latex](0,0)--(1.5,1.5);
		\draw[-latex](0,0)--(1.5,-1.5);
		\draw[-latex](-3,3)--(-1.5,1.5);
		
		\coordinate (I) at (0.2,0);
		
		\coordinate (E) at (0,1);
		\fill (E) circle (0pt) node[above] {\large $\Omega_2$};
		\coordinate (D) at (1.2,0.6);
		\fill (D) circle (0pt) node[right] {\large $\Omega_1$};
		\coordinate (F) at (0,-1);
		\fill (F) circle (0pt) node[below] {\large $\Omega_5$};
		\coordinate (J) at (-1.2,-0.6);
		\fill (J) circle (0pt) node[left] {\large $\Omega_4$};
		\coordinate (k) at (-1.2,0.6);
		\fill (k) circle (0pt) node[left] {\large $\Omega_3$};
		\coordinate (k) at (1.2,-0.6);
		\fill (k) circle (0pt) node[right] {\large$\Omega_6$};
		\fill (I) circle (0pt) node[below] {$-\xi$};

		\coordinate (A) at (0,2);
		\coordinate (B) at (-2,2);
		\coordinate (C) at (1.5,-1.8);
		\coordinate (D) at (-3.5,-1.8);
		\coordinate (E) at (-1,1);
		\coordinate (F) at (-1,-2);
		\draw[red] (E) circle (0.12);
		\draw[red] (F) circle (0.12);	
		\draw[red] (A) circle (0.12);
		\draw[red] (B) circle (0.12);	
		\draw[red] (C) circle (0.12);		
		\draw[red] (D) circle (0.12);
		
		\draw[blue] (0,0) circle (0.3);
		
		\fill (A) circle (1pt) node[right] {$\omega_n$};
		\fill (B) circle (1pt) node[left] {$-\overline{\omega_n}$};
		\fill (C) circle (1pt) node[right] {$\gamma_n$};
		\fill (D) circle (1pt) node[left] {$-\overline{\gamma_n}$};
		\fill (E) circle (1pt) node[right] {$\eta_m$};
		\fill (F) circle (1pt) node[right] {$\tau_m$};
		
	\end{tikzpicture}
	\caption{\small The $ \mathcal{X}(k)$ and  $\mathcal{Y}(k)$ are supported in the red and blue circle respectively. }
	\label{figR2}
\end{figure}

On each $\Omega_j$,  we define smooth functions $ R_j(k)$  as follows
\begin{align*}
	R_1(k)=(1-\mathcal{X}(k))&\left(r_1({\rm Re}k)\delta^{-2}(k)+\sin(2\phi)\left(f_1(k)-r_1({\rm Re}k)\delta^{-2}(k)\right)\right),\\
	R_3(k)=(1-\mathcal{X}(k))&\left[ \frac{\sigma r_2}{1+\sigma r_1 r_2}({\rm Re}k)\delta^{2}(k)\left(1-\sin(2\phi)\right)+ \sin(2\phi)f_3(k)\right],\\
	R_4(k)=(1-\mathcal{X}(k))&\left[ \frac{ r_1}{1+\sigma r_1 r_2}({\rm Re}k)\delta^{-2}(k)\left(1-\sin(2\phi)\right)+ \sin(2\phi)f_4(k)\right],\\
	R_6(k)=(1-\mathcal{X}(k))&\left(\sigma r_2({\rm Re}k)\delta^{2}(k)+\sin(2\phi)\left(f_6(k)-\sigma r_2({\rm Re}k)\delta^{2}(k)\right)\right),
\end{align*}
where
\begin{align*}
	f_1(k)&=\mathcal{Y}(k)r_1(-\xi)\delta_0^{-2}(\xi+k)^{-2{\rm i}\nu(-\xi)},\\
	f_3(k)&=\mathcal{Y}(k)\frac{\sigma r_2}{1+\sigma r_1 r_2}(-\xi)\delta_0^{2}(\xi+k)^{2{\rm i}\nu(-\xi)},\\
	f_4(k)&=\mathcal{Y}(k)\frac{r_1}{1+\sigma r_1 r_2}(-\xi)\delta_0^{-2}(\xi+k)^{-2{\rm i}\nu(-\xi)},\\
	f_6(k)&=\mathcal{Y}(k)\sigma r_2(-\xi)\delta_0^{2}(\xi+k)^{2{\rm i}\nu(-\xi)},
\end{align*}
and\begin{equation}
	R_2(k)= R_5(k)=I.\nonumber
\end{equation}
It is obvious that $R_j$ has  boundary value  with
\begin{align*}
	R_1(k)&=\left\{ \begin{array}{ll}
		f_1(k),&k\in \Sigma_1,\\[12pt]
		r_1({\rm Re}k)\delta^{-2}(k),&k>-\xi,\\
	\end{array}\right.\\
	R_3(k)&=\left\{ \begin{array}{ll}
		f_3(k),&k\in \Sigma_3,\\[12pt]
		\frac{\sigma r_2}{1+\sigma r_1 r_2}({\rm Re}k)\delta^{2}(k),&k<-\xi,\\
	\end{array}\right.\\
	R_4(k)&=\left\{ \begin{array}{ll}
		f_4(k),&k\in \Sigma_4,\\[12pt]
		\frac{ r_1}{1+\sigma r_1 r_2}({\rm Re}k)\delta^{-2}(k),&k<-\xi,\\
	\end{array}\right.\\
	R_6(k)&=\left\{ \begin{array}{ll}
		f_6(k),&k\in \Sigma_6,\\[12pt]
		\sigma r_2({\rm Re}k)\delta^{2}(k),&k>-\xi.\\
	\end{array}\right.\\
\end{align*}
\begin{proposition} \label{dbar-R2}
	$ R_j(k)$ satisfies the following properties
	\begin{align}
		|\overline{\partial}R_1(k)|&\lesssim|\overline{\partial}\mathcal{X}|+|\overline{\partial}\mathcal{Y}|+|r'_1({\rm Re}k)|+|k+\xi|^{-\alpha},\nonumber\\
		|\overline{\partial}R_3(k)|&\lesssim|\overline{\partial}\mathcal{X}|+|\overline{\partial}\mathcal{Y}|+\left|\left( \frac{\sigma r_2}{1+\sigma r_1 r_2}\right) '({\rm Re}k)\right|+|k+\xi|^{-\alpha},\nonumber\\
		|\overline{\partial}R_4(k)|&\lesssim|\overline{\partial}\mathcal{X}|+|\overline{\partial}\mathcal{Y}|+\left|\left( \frac{r_1}{1+\sigma r_1 r_2}\right) '({\rm Re}k)\right|+|k+\xi|^{-\alpha},\nonumber\\
		|\overline{\partial}R_6(k)|&\lesssim|\overline{\partial}\mathcal{X}|+|\overline{\partial}\mathcal{Y}|+|r'_2({\rm Re}k)|+|k+\xi|^{-\alpha},\nonumber
	\end{align}
	where
	\begin{align}
		\alpha=\left\lbrace\begin{array}{ll}
			\frac{1}{2}+{\rm Im}\nu(-\xi)&,{\rm Im}\nu(-\xi)\in\left[0,\frac{1}{2}\right),\\[12pt]
			\frac{1}{2}&,{\rm Im}\nu(-\xi)\in\left(-\frac{1}{2},0\right).
		\end{array}\right.\label{alpha}
	\end{align}
\end{proposition}
\begin{proof}
	We only give  details of the proof for $R_1 $. The others can be demonstrated in the same way.
	\begin{align}
		|\overline{\partial}R_1(k)|=&-\overline{\partial}\mathcal{X}R_1(k)+(1-\mathcal{X})\frac{1}{2}r'_1({\rm Re}k)\delta^{-2}(k)(1-\sin(2\phi))\nonumber\\
		&+(1-\mathcal{X})\overline{\partial}\mathcal{Y}r_1(-\xi)\delta_0^{-2}(\xi+k)^{-2{\rm i}\nu(-\xi)}\sin(2\phi)\nonumber\\
		&+\frac{i{\rm e}^{i\phi}}{\rho}(1-\mathcal{X})\cos(2\phi)\left(f_1(k)-r_1({\rm Re}k)\delta^{-2}(k)\right) \nonumber\\
		\lesssim&|\overline{\partial}\mathcal{X}|+|\overline{\partial}\mathcal{Y}|+|r'_1({\rm Re}k)|+|\xi+k|^{-1}r_1({\rm Re}k)\delta^{-2}(k)\chi_{\{\rho>\rho_0\}}\nonumber\\
		&+|\xi+k|^{-1}\left(r_1(-\xi)\delta_0^{-2}(\xi+k)^{-2{\rm i}\nu(-\xi)}-r_1({\rm Re}k)\delta^{-2}(k)\right)\chi_{\{\rho<\frac{3}{2}\rho_0\}} \nonumber\\
		\lesssim&|\overline{\partial}\mathcal{X}|+|\overline{\partial}\mathcal{Y}|+|r'_1({\rm Re}k)|+|\xi+k|^{-1}\chi_{\{|\rho|>\rho_0\}}\nonumber\\
		&+\left(|k+\xi|^{-\frac{1}{2}-{\rm Im}\nu(-\xi)}+|k+\xi|^{-\frac{1}{2}}\right)\chi_{\{|\rho|<\frac{3}{2}\rho_0\}}.\nonumber
	\end{align}
	Depending on (\ref{1/2}) and $r_1\in H^1(\mathbb{R})$, then $r_1 $ is $ \frac{1}{2}$-H\"{o}lder continuous. So
	\begin{align*}
		&\left(r_1(-\xi)\delta_0^{-2}(\xi+k)^{-2{\rm i}\nu(-\xi)}-r_1({\rm Re}k)\delta^{-2}(k)\right)\chi_{\{\rho<\frac{3}{2}\rho_0\}}\\
		\leq&\left(|r_1(-\xi)||\delta_0^{-2}(\xi+k)^{-2{\rm i}\nu(-\xi)}-\delta^{-2}(k)|+|r_1(\xi)-r_1({\rm Re}k)|\delta^{-2}(k)\right)\chi_{\{\rho<\frac{3}{2}\rho_0\}}\\
		\lesssim&\left(|k+\xi|^{\frac{1}{2}-{\rm Im}\nu(-\xi)}+|k+\xi|^{\frac{1}{2}}\right)\chi_{\{\rho<\frac{3}{2}\rho_0\}}.
	\end{align*}
	Meanwhile, $\overline{\partial}\mathcal{Y}$ has bounded support on $\{\rho_0<|\rho|<\frac{3}{2}\rho_0\}$, so that
	\begin{align*}
		|\overline{\partial}\mathcal{Y}||k+\xi|^{-2|{\rm Im}\nu(-\xi)|}\lesssim|\overline{\partial}\mathcal{Y}|.
	\end{align*}
	Under the assumption $|{\rm Im}\nu(-\xi)|<\frac{1}{2}$,
	\begin{equation}
		|\xi+k|^{-1}\chi_{\{|\rho|>\rho_0\}}\lesssim|k+\xi|^{-\frac{1}{2}-{\rm Im}\nu(-\xi)}+|k+\xi|^{-\frac{1}{2}}.\nonumber
	\end{equation}
	Therefore,
	\begin{align}
		|\overline{\partial}R_1(k)|
		\lesssim|&\overline{\partial}\mathcal{X}|+|\overline{\partial}\mathcal{Y}|+|r'_1({\rm Re}k)|+|\xi+k|^{-\alpha}.\nonumber
	\end{align}
\end{proof}
To make continuous extension for  the jump matrix $\tilde{J}$  to remove the jump from $\Sigma$, let
\begin{equation}
	\mathcal{R}(z)=\left\{\begin{array}{lll}
		\left(\begin{array}{cc}
			1 & (-1)^jR_j(z){\rm e}^{-2{\rm i}t\theta}\\
			0 & 1
		\end{array}\right), & z\in \Omega_j,j=3,6;\\
		\\
		\left(\begin{array}{cc}
			1 & 0\\
			(-1)^jR_j(z){\rm e}^{2{\rm i}t\theta} & 1
		\end{array}\right),  &z\in \Omega_j,j=1,4;\\
		\\
		I,  &elsewhere,\\
	\end{array}\right.\label{R(2)}
\end{equation}
Besides, the new problem is hoped to  take advantage of the decay/growth of ${\rm e}^{2{\rm i}t\theta(z)}$ for $z\notin\Sigma$. We give the second transform
\begin{align}
	\tilde{M}(x,t;k)\mathcal{R}(k)=M^{(2)}(x,t;k).
\end{align}
Because $\mathcal{R}$ is a sectionally continuous function,  $M^{(2)}(x,t;k) $ is derivable on $\mathbb{C}\setminus\Sigma$ with $\Sigma=\cup_{j=1}^4\Sigma_j$. We can derive the derivative condition
\begin{equation}
	\overline{\partial}M^{(2)}=\overline{\partial}(\tilde{M}\mathcal{R})=\tilde{M}\overline{\partial}\mathcal{R}=M^{(2)}\overline{\partial}\mathcal{R}.
\end{equation}
Notice that  $\mathcal{R}=I $ near the pole of $\tilde{M} $. The matrix valued function $M^{(2)}(x,t;k) $ satisfies following mixed $ \overline{\partial}$-RH problem

\noindent\textbf{RHP2}.  Find a matrix-valued function $M^{(2)}(x,t;k)$ which satisfies:

$\blacktriangleright$ $M^{(2)}(x,t;k)$ is continuous in $\mathbb{C}\ \setminus\Sigma $ and  meromorphic in $\Omega_2\cup\Omega_5$.

$\blacktriangleright$  $M^{(2)}(x,t;k)$ has continuous boundary values $M^{(2)}(x,t;k^\pm)$ on $\Sigma$ with
\begin{align}
	M^{(2)}(x,t;k^+)&=M^{(2)}(x,t;k^-)J^{(2)}(x,t;k),\end{align}
where
\begin{align}
	J^{(2)}(x,t;k)&=\left\{ \begin{array}{llll}
		\left(\begin{array}{cc}
			1& 0\ \\
			f_1(k){\rm e}^{2{\rm i}t\theta}&1\\
		\end{array}\right),  & k\in\Sigma_1,\\[12pt]
		\left(\begin{array}{cc}
			1& f_3(k){\rm e}^{-2{\rm i}t\theta} \\
			0 &1\\
		\end{array}\right),&k\in\Sigma_2, \\
		\left(\begin{array}{cc}
			1&0 \\
			f_4(k){\rm e}^{2{\rm i}t\theta} &1\\
		\end{array}\right),&k\in\Sigma_3,\\
		\left(\begin{array}{cc}
			1&f_6(k){\rm e}^{-2{\rm i}t\theta}\\
			0  &1\\
		\end{array}\right) , &k\in\Sigma_4.\\
	\end{array}\right.
\end{align}

$\blacktriangleright$   The residue conditions of $ M^{(2)}(x,t;k)$ satisfy (\ref{res1+}), (\ref{res1-}) by  replacing $ \tilde{M}(x,t;k)$.

$\blacktriangleright$  $\bar\partial$-derivative condition
\begin{equation}
	\overline{\partial}M^{(2)}=M^{(2)}\overline{\partial}\mathcal{R}.\label{dbar-m2}
\end{equation}

$\blacktriangleright$  Asymptotic condition: $M^{(2)}(x,t;k)=I+\mathcal{O}(k^{-1}),\  k\rightarrow \infty.$

To solve the mixed RH problem, we define a function $M^{RHP}(x,t;k) $ satisfies following model Riemann-Hilbert problem   with $\bar\partial \mathcal{R}\equiv0$:

\noindent\textbf{RHP3}.  Find a matrix-valued function $M^{RHP}(x,t;k)$ which satisfies:

$\blacktriangleright$ $M^{RHP}(x,t;k) $ is meromorphic in $\mathbb{C}\ \setminus \Sigma$.

$\blacktriangleright$   $M^{RHP}(x,t;k)$ has continuous boundary values $M_{RHP}(x,t;k^\pm)$ on $\Sigma$ with
\begin{align}
	M^{RHP}_+(x,t;k )&= M^{RHP}_-(x,t;k )J^{(2)}(x,t;k),\end{align}
where
\begin{align}
	J^{(2)}(x,t;k)&=\left\{ \begin{array}{llll}
		\left(\begin{array}{cc}
			1& 0\ \\
			f_1(k){\rm e}^{2{\rm i}t\theta}&1\\
		\end{array}\right),  & k\in\Sigma_1,\\[12pt]
		\left(\begin{array}{cc}
			1& f_3(k){\rm e}^{-2{\rm i}t\theta} \\
			0 &1\\
		\end{array}\right),&k\in\Sigma_2, \\
		\left(\begin{array}{cc}
			1&0 \\
			f_4(k){\rm e}^{2{\rm i}t\theta} &1\\
		\end{array}\right),&k\in\Sigma_3,\\
		\left(\begin{array}{cc}
			1&f_6(k){\rm e}^{-2{\rm i}t\theta}\\
			0  &1\\
		\end{array}\right) , &k\in\Sigma_4.\\
	\end{array}\right.
\end{align}

$\blacktriangleright$ Residue condition: The residue conditions of $M^{RHP}(x,t;k)$ satisfy (\ref{res1+}), (\ref{res1-}) by  replacing $ \tilde{M}(x,t;k)$.

$\blacktriangleright$ Asymptotic condition:
\begin{equation}
	M^{RHP}(x,t;k)=I+\mathcal{O}(k^{-1}),\  k\rightarrow \infty.
\end{equation}

To prove the existence of function $ M^{RHP}(x,t;k)$, we divide it to two parts:
\begin{align}
	M^{RHP}(x,t;k)&=\left\{ \begin{array}{ll}
		E(x,t;k)M_{sol}(x,t;k)M^{PC}(z),  & |k+\xi|<\rho_0,\\[12pt]
		E(x,t;k)M_{sol}(x,t;k),&|k+\xi|>\rho_0, \\
	\end{array}\right.
\end{align}
where  the matrix function $ M_{sol}(x,t;k)$ satisfies (\ref{res1+}), (\ref{res1-}) with $ M_{sol}(x,t;k)$ replacing $ \tilde{M}(x,t;k)$. And analytical  in the elsewhere of $ k\in\mathbb{C}$. And\begin{equation}
	M_{sol}=I+\mathcal{O}(k^{-1}),\  k\rightarrow \infty.
\end{equation}
which will be discussed in Section \ref{msol}.

While  $ M^{PC}$ is the well known parabolic cylinder model satisfying  jump condition \cite{PC1}
\begin{align}
	J^{PC}(x,t;k)&=\left\{ \begin{array}{llll}
		\left(\begin{array}{cc}
			1& 0\ \\
			r_{\xi}z^{-2{\rm i}\nu(-\xi)}{\rm e}^{\frac{{\rm i}}{2}z^2}&1\\
		\end{array}\right),  & z\in\mathbb{C}, \arg(z)=\frac{\pi}{4},\\[12pt]
		\left(\begin{array}{cc}
			1& \frac{\check{r}_{\xi}}{1+r_{\xi}\check{r}_{\xi}}z^{2{\rm i}\nu(-\xi)}{\rm e}^{-\frac{{\rm i}}{2}z^2}\\
			0 &1\\
		\end{array}\right),& z\in\mathbb{C}, \arg(z)=\frac{3\pi}{4}, \\
		\left(\begin{array}{cc}
			1&0 \\
			\frac{r_{\xi}}{1+r_{\xi}\check{r}_{\xi}}z^{-2{\rm i}\nu(-\xi)}{\rm e}^{\frac{{\rm i}}{2}z^2} &1\\
		\end{array}\right),& z\in\mathbb{C}, \arg(z)=\frac{5\pi}{4},\\
		\left(\begin{array}{cc}
			1&\check{r}_{\xi}z^{2{\rm i}\nu(-\xi)}{\rm e}^{-\frac{{\rm i}}{2}z^2}\\
			0  &1\\
		\end{array}\right) , & z\in\mathbb{C}, \arg(z)=\frac{7\pi}{4},\\
	\end{array}\right.
\end{align}
with parameters
\begin{align}
	r_{\xi}=&r_1(-\xi)\delta_0^{-2}(8t)^{i\nu(-\xi)}{\rm e}^{-4it\xi^2},\
	\check{r}_{\xi}=\sigma r_2(-\xi)\delta_0^{2}(8t)^{-i\nu(-\xi)}{\rm e}^{4it\xi^2} .
\end{align}

In order to match   the model, let $z = z(k)$ denote the resealed local variable
\begin{align}
	z=&\sqrt{8t}(\xi+k).
\end{align}
\begin{proposition}The large-$z$ asymptotics  of $M^{PC}(x,t;z)$  satisfies \label{pcasy}
	\begin{align}
		M^{PC}(x,t;z)&=I+\frac{1}{z}\left(\begin{array}{cc}
			0	&-{\rm i}\beta_{12}\\
			{\rm i}\beta_{21}&0\\
		\end{array}\right)+\mathcal{O}(z^{-2})\nonumber\\
		\beta_{12}&=\frac{\sqrt{2\pi}{\rm e}^{-\frac{\pi}{2}\nu(-\xi)}{\rm e}^{\frac{\pi{\rm i}}{4}}}{r_{\xi}\Gamma(-{\rm i}\nu(-\xi))}\triangleq t^{{\rm Im}\nu(-\xi)} \tilde{\beta}_{12}\label{b12},\\
		\beta_{21}&=-\frac{\sqrt{2\pi}{\rm e}^{-\frac{\pi}{2}\nu(-\xi)}{\rm e}^{-\frac{\pi{\rm i}}{4}}}{\check{r}_{\xi}\Gamma({\rm i}\nu(-\xi))}\triangleq t^{-{\rm Im}\nu(-\xi)} \tilde{\beta}_{21}.
	\end{align}
\end{proposition}
\begin{proof}
	The proof of the proposition  was given in  \cite{Liu3}.
\end{proof}
Because $ \mathcal{Y}(k)=1$ when $ |\xi+k|<\rho_0$, we have
\begin{align}
	J^{(2)}(x,t;k)&=\left\{ \begin{array}{llll}
		\left(\begin{array}{cc}
			1& 0\ \\
			r_{\xi}z^{-2{\rm i}\nu(-\xi)}{\rm e}^{\frac{{\rm i}}{2}z^2}&1\\
		\end{array}\right),  & |k+\xi|<\rho_0,k\in\Sigma_1,\\[12pt]
		\left(\begin{array}{cc}
			1& \frac{\check{r}_{\xi}}{1+r_{\xi}\check{r}_{\xi}}z^{2{\rm i}\nu(-\xi)}{\rm e}^{-\frac{{\rm i}}{2}z^2}\\
			0 &1\\
		\end{array}\right),& |k+\xi|<\rho_0,k\in\Sigma_2, \\
		\left(\begin{array}{cc}
			1&0 \\
			\frac{r_{\xi}}{1+r_{\xi}\check{r}_{\xi}}z^{-2{\rm i}\nu(-\xi)}{\rm e}^{\frac{{\rm i}}{2}z^2} &1\\
		\end{array}\right),& |k+\xi|<\rho_0,k\in\Sigma_3,\\
		\left(\begin{array}{cc}
			1&\check{r}_{\xi}z^{2{\rm i}\nu(-\xi)}{\rm e}^{-\frac{{\rm i}}{2}z^2}\\
			0  &1\\
		\end{array}\right) , & |k+\xi|<\rho_0,k\in\Sigma_4.\\
	\end{array}\right.
\end{align}
Comparing with the jump condition of $M^{PC}(z)$,  the error function $ E(x,t;k)$ satisfies  a small norm RH problem

\noindent\textbf{RHP4}.  Find a matrix-valued function $E(x,t;k)$ which satisfies

$\blacktriangleright$ $E(x,t;k)$ is analytic in $\mathbb{C}\ \setminus\Sigma^{E}, $
where as shown in Figure \ref{figE}, \begin{equation}
	\widehat{\Sigma}_j=\Sigma_j\setminus\left\{k\in\mathbb{C}\big||k+\xi|<\rho_0\right\},\  \Sigma^{E}=\bigcup_{j=1}^4\widehat{\Sigma}_j\cup\left\{k\in\mathbb{C}\big||k+\xi|=\rho_0\right\}.
\end{equation}

$\blacktriangleright$ Jump condition:
For $k\in\Sigma^{E},$ \begin{equation}
	E_+(x,t;k )=E_-(x,t;k  )J^{E}(x,t;k),
\end{equation}where
\begin{align}
	J^{E}(x,t;k)&=\left\{ \begin{array}{ll}
		M_{sol}(x,t;k)J^{(2)}(x,t;k)M_{sol}^{-1}(x,t;k),&k\in\widehat{\Sigma}_j,\\
		M_{sol}(x,t;k)M^{PC}(x,t;k)M_{sol}^{-1}(x,t;k), &k\in\left\{k\in\mathbb{C}\big||k+\xi|=\rho_0\right\}.\\
	\end{array}\right. \nonumber
\end{align}

$\blacktriangleright$ Asymptotic condition:
\begin{equation}
	E(x,t;k)=I+\mathcal{O}(k^{-1}),\  k\rightarrow \infty.
\end{equation}

\begin{figure}[H]
	\centering
	\begin{tikzpicture}[scale=0.7]
		\draw(1,1)--(3,3)node[black][right]{$\widehat{\Sigma}_1$};
		\draw(-1,1)--(-3,3)node[left]{$\widehat{\Sigma}_2$};
		\draw(-1,-1)--(-3,-3)node[left]{$\widehat{\Sigma}_3$};
		\draw(1,-1)--(3,-3)node[right]{$\widehat{\Sigma}_4$};
		\draw[dashed](-4,0)--(4,0)node[right]{ Re$k$};
		\draw[dashed](0,-3.5)--(0,3.5)node[above]{ $-\xi$};
		\draw[-latex](-3,-3)--(-1.5,-1.5);
		\draw[-latex](1,1)--(1.5,1.5);
		\draw[-latex](1,-1)--(1.5,-1.5);
		\draw[-latex](-3,3)--(-1.5,1.5);
		\draw (0,0) circle (1.43) ;
		\draw[-latex] (1.43,0) arc (0:120:1.43);
		\fill (1.5,1) circle(0pt) node[right]{$\Sigma^{E}$};
	\end{tikzpicture}
	\caption{ The jump contour of $E(x,t;k)$}\label{figE}
\end{figure}
To prove the existence of $M_{RHP}(x,t;k) $, we only need to find $E(x,t;k)$ for RHP4.  $\mathcal{C}_{\pm} $ are the  limit of general Cauchy operators:
\begin{equation}
	\mathcal{C}_{\pm}(f)(s)=\lim_{z\to \Sigma^{E}_\pm}\frac{1}{2\pi{\rm i}}\int_{\Sigma^{E}}\dfrac{f(s)}{s-z}ds.
\end{equation}
Let
\begin{equation}
	\mathcal{C}_{\omega}[f]=\mathcal{C}_{-}(f(J^{E}-I)).
\end{equation}
The RHP4 is solvable if and only if there is a function $\mu $ satisfies
\begin{equation}
	(Id-\mathcal{C}_{\omega})(I+\mu)=I.
\end{equation}
\begin{proposition}
	$\mathcal{C}_{\omega}$ is a bounded operator on $L^{2}(\Sigma^{E})\to L^{2}(\Sigma^{E})$ with:
	\begin{equation}
		||\mathcal{C}_{\omega}||_{\mathcal{B}(L^{2}(\Sigma^{E}))}=
		\mathcal{O}(t^{-\frac{1}{2}+|{\rm Im}\nu(-\xi)|}).
	\end{equation}
\end{proposition}
\begin{proof}
	\begin{align}
		||\mathcal{C}_{\omega}[f]||_{L^{2}(\Sigma^{E})}
		&\lesssim||f(J^{E}-I)||_{L^{2}(\Sigma^{E})}\nonumber\\
		&\leq||f||_{L^{2}(\Sigma^{E})}||J^{E}-I||_{L^{\infty}(\Sigma^{E})},\nonumber\\
		||\mathcal{C}_{\omega}||_{\mathcal{B}(L^{2}(\Sigma^{E}))}
		&\lesssim||J^{E}-I||_{L^{\infty}(\Sigma^{E})}\nonumber\\
		&\lesssim||J^{E}-I||_{L^{\infty}(\bigcup_{j=1}^{4}\Sigma_j)}+||J^{E}-I||_{L^{\infty}(\left\lbrace |k+\xi|=\rho_0\right\rbrace )},\nonumber\\
		&\leq\mathcal{O}(t^{|{\rm Im}\nu(-\xi)|}\exp(-4\rho_0^2t))+\mathcal{O}(t^{-\frac{1}{2}+|{\rm Im}\nu(-\xi)|})\nonumber\\
		&=\mathcal{O}(t^{-\frac{1}{2}+|{\rm Im}\nu(-\xi)|}).\nonumber
	\end{align}
\end{proof}
Then for sufficiently large $t$, $Id-\mathcal{C}_{\omega} $ becomes a bijection in $L^2(\Sigma^{E})$.  Specially, equation
\begin{equation}
	(Id-\mathcal{C}_{\omega})\mu=\mathcal{C}_{\omega}I,
\end{equation} has a solution $\mu\in L^2(\Sigma^{E})$
and satisfying \begin{equation}
	||\mu||_{L^2(\Sigma^{E})}\lesssim||\mathcal{C}_{\omega}I||_{L^2(\Sigma^{E})}\lesssim\mathcal{O}(t^{-\frac{1}{2}+|{\rm Im}\nu(-\xi)|}).
\end{equation}
RHP4 can be solved by
\begin{equation}
	E(x,t;k)=I+\frac{1}{2\pi{\rm i}}\int_{\Sigma^E}\frac{(I+\mu)(J^E-I)(s)}{s-k}ds,
\end{equation}
\begin{proposition}
	$E(x,t;k)$ satisfies large $k$ asymptotic condition
	\begin{align}
		E(x,t;k)&=I+\frac{E^{(1)}(x,t)}{k}+\mathcal{O}(k^{-2}),\\
		E^{(1)}(x,t)&=\frac{1}{\sqrt{8t}}\left(\begin{array}{cc}
			0&-{\rm i}t^{{\rm Im}\nu(-\xi)}\tilde{\beta}_{l2}\\
			{\rm i}t^{-{\rm Im}\nu(-\xi)}\tilde{\beta}_{21}&0\\
		\end{array}\right)\nonumber\\
		&+\mathcal{O}(t^{-\frac{1}{2}+|{\rm Im}\nu(-\xi)|})\begin{array}{cc}
			\Big(\mathcal{O}(t^{-\frac{1}{2}-{\rm Im}\nu(-\xi)}),\mathcal{O}(t^{-\frac{1}{2}+{\rm Im}\nu(-\xi)})\Big)\end{array}.\label{Easy}
	\end{align}
\end{proposition}
\begin{proof}According to Proposition \ref{pcasy} and Cauchy integral formula, we deduce
	\begin{align*}
		E^{(1)}(x,t)&=\frac{1}{2\pi{\rm i}}\int_{\Sigma^E}(I+\mu)(J^E-I)(s)ds\\
		&=\frac{1}{2\pi{\rm i}}\int_{\Sigma^E}(J^E-I)(s)ds+\frac{1}{2\pi{\rm i}}\int_{\Sigma^E}\mu(J^E-I)(s)ds\\
		&=\frac{1}{\sqrt{8t}}\left(\begin{array}{cc}
			0&-{\rm i}t^{{\rm Im}\nu(-\xi)}\tilde{\beta}_{l2}\\
			{\rm i}t^{-{\rm Im}\nu(-\xi)}\tilde{\beta}_{21}&0\\
		\end{array}\right)+\mathcal{O}(\exp(-4\rho_0^2t))\\
		&+\mathcal{O}(t^{-\frac{1}{2}+|{\rm Im}\nu(-\xi)|})\begin{array}{cc}
			\Big(\mathcal{O}(t^{-\frac{1}{2}-{\rm Im}\nu(-\xi)}),\mathcal{O}(t^{-\frac{1}{2}+{\rm Im}\nu(-\xi)})\Big),\end{array}		
	\end{align*}
	where
	\begin{align*}
		&\frac{1}{2\pi{\rm i}}\int_{\Sigma^E}\mu(J^E-I)(s)ds\leq
		||\mu||_{L^2(\Sigma^{E})}||J^E-I||_{L^2(\Sigma^{E})}\\
		&=\mathcal{O}(t^{-\frac{1}{2}+|{\rm Im}\nu(-\xi)|})\begin{array}{cc}
			\Big(\mathcal{O}(t^{-\frac{1}{2}-{\rm Im}\nu(-\xi)}),\mathcal{O}(t^{-\frac{1}{2}+{\rm Im}\nu(-\xi)})\Big).		
		\end{array}
	\end{align*}
\end{proof}

\subsection{Analysis on a  soliton  model    } \label{msol}
\quad
Recall the definition of $M_{sol}$ as follow.

\noindent\textbf{RHP5}.  Find a matrix-valued function $M_{sol}(x,t;k)$ which satisfies

$\blacktriangleright$ $M_{sol}(x,t;k) $ is meromorphic in $\mathbb{C}$.

$\blacktriangleright$ Asymptotic condition:
\begin{equation}
	M_{sol}(x,t;k)=I+\mathcal{O}(k^{-1}),\  k\rightarrow \infty.
\end{equation}

$\blacktriangleright$ Residue condition: \begin{align}
	\begin{aligned}
		\res_{k=\omega_{n} }M_{sol}(x,t;k)&=\lim_{k\rightarrow \omega_n}M_{sol}(x,t;k)\left(\begin{array}{cc}
			0 & 0 \\
			\delta^{-2}(\omega_{n})c_n{\rm e}^{2{\rm i}t\theta} & 0
		\end{array}\right),\\
		\res_{k=\eta_{m} }M_{sol}(x,t;k)&=\lim_{k\rightarrow \eta_{m}}M_{sol}(x,t;k)\left(\begin{array}{cc}
			0 & 0 \\
			\delta^{-2}(\eta_{m})e_m{\rm e}^{2{\rm i}t\theta} & 0
		\end{array}\right),
	\end{aligned}\\
	\begin{aligned}
		\res_{k=\gamma_{n} }M_{sol}(x,t;k)&=\lim_{k\rightarrow \gamma_{n}}M_{sol}(x,t;k)\left(\begin{array}{cc}
			0 & \delta^2(\gamma_{n}) d_n{\rm e}^{-2{\rm i}t\theta} \\
			0 & 0
		\end{array}\right),\\
		\res_{k=\tau_{m} }M_{sol}(x,t;k)&=\lim_{k\rightarrow \tau_{m}}M_{sol}(x,t;k)\left(\begin{array}{cc}
			0 & \delta^2(\tau_{m}) f_m{\rm e}^{-2{\rm i}t\theta} \\
			0 & 0
		\end{array}\right).
	\end{aligned}
\end{align}

Let  $\alpha_\omega$, $\beta_\gamma$ denote the column valued residue of $M_{sol}$ whose subscripts indicate the pole. Then $M_{sol}$ can be written as
\begin{equation}
	M_{sol}(x,t;k)=I+\sum_{\omega\in \mathcal{N}_1}\frac
	{1}{k-\omega}
	\left(\begin{array}{cc}
		\alpha_{\omega}^{(1)}&0\\
		\alpha_{\omega}^{(2)}&0\\
	\end{array}\right)
	+\sum_{\gamma\in \mathcal{N}_2}\frac
	{1}{k-\gamma}
	\left(\begin{array}{cc}
		0&\beta_\gamma^{(1)}\\
		0&\beta_\gamma^{(2)}\\
	\end{array}\right).
\end{equation}
In above expression, $\alpha_\omega,\ \omega\in\mathcal{N}_1$ and $\beta_\gamma,\ \gamma\in\mathcal{N}_2$ are the solution of the following linear equations:
\begin{align}\label{res1}
	\begin{aligned}
		\left[I
		+\sum_{\gamma\in \mathcal{N}_2}\frac
		{1}{\omega_n-\gamma}
		\left(\begin{array}{cc}
			0&\beta_\gamma^{(1)}\\
			0&\beta_\gamma^{(2)}\\
		\end{array}\right)\right]
		\left(\begin{array}{cc}
			0 &0  \\
			c_n\delta(\omega_n)^{-2}{\rm e}^{2{\rm i}t\theta} & 0
		\end{array}\right)=&
		\left(\begin{array}{cc}
			\alpha_{\omega_{n}}^{(1)}&0\\
			\alpha_{\omega_{n}}^{(2)}&0\\
		\end{array}\right),n=\pm1,...,\pm N_1,\\
		\left[I
		+\sum_{\gamma\in \mathcal{N}_2}\frac
		{1}{\eta_m-\gamma}
		\left(\begin{array}{cc}
			0&\beta_\gamma^{(1)}\\
			0&\beta_\gamma^{(2)}\\
		\end{array}\right)\right]
		\left(\begin{array}{cc}
			0 &0  \\
			e_m\delta^{-2}(\eta_{m}){\rm e}^{2{\rm i}t\theta} & 0
		\end{array}\right)=&
		\left(\begin{array}{cc}
			\alpha_{\eta_{m}}^{(1)}&0\\
			\alpha_{\eta_{m}}^{(2)}&0\\
		\end{array}\right),m=1,...,M_1,
	\end{aligned}\\\label{res2}
	\begin{aligned}
		\left[I+\sum_{\omega\in \mathcal{N}_1}\frac
		{1}{\gamma_n-\omega}
		\left(\begin{array}{cc}
			\alpha_{\omega}^{(1)}&0\\
			\alpha_{\omega}^{(2)}&0\\
		\end{array}\right)\right]
		\left(\begin{array}{cc}
			0 &d_n\delta(\gamma_n)^2{\rm e}^{-2{\rm i}t\theta}  \\
			0& 0
		\end{array}\right)=&
		\left(\begin{array}{cc}
			0&\beta_{\gamma_n}^{(1)}\\
			0&\beta_{\gamma_n}^{(2)}\\
		\end{array}\right),n=\pm1,...,\pm N_2.\\
		\left[I
		+\sum_{\omega\in \mathcal{N}_1}\frac
		{1}{\tau_m-\omega}
		\left(\begin{array}{cc}
			\alpha_{\omega}^{(1)}&0\\
			\alpha_{\omega}^{(2)}&0\\
		\end{array}\right)\right]
		\left(\begin{array}{cc}
			0 &f_m\delta^{2}(\tau_{m}){\rm e}^{-2{\rm i}t\theta} \\
			0  & 0
		\end{array}\right)=&
		\left(\begin{array}{cc}
			0&\beta_{\tau_m}^{(1)}\\
			0&\beta_{\tau_m}^{(2)}\\
		\end{array}\right),m=1,...,M_2,
	\end{aligned}
\end{align}
The existence of $ M_{sol}(x,t;k)$ can be seen by the coefficients matrix  of above linear algebraic equation is reversible.
As a additional condition of the Riemann-Hilbert problem approaching to the NNLS equation, the symmetry (\ref{sym}) implies
\begin{align}
	M(x,t;k)=\left\{ \begin{array}{ll}
		\Lambda\overline{M(-x,t;-\bar{k})}\Lambda^{-1}\left(\begin{array}{ll}
			\frac{1}{a_1(k)}&0\\
			0&a_1(k)
		\end{array}\right) , &k\in\mathbb{C}^+;\\
		\Lambda\overline{M(-x,t;-\bar{k})}\Lambda^{-1}\left(\begin{array}{ll}
			a_2(k)&0\\
			0&\frac{1}{a_2(k)}
		\end{array}\right) , &k\in\mathbb{C}^-.\\
	\end{array}\right.\label{Msym}
\end{align}
The $a_1,\ a_2$ here are given by the scattering data $r_1,\  r_2$ via the well-known trace formula
\begin{align}
	a_1(k)=\frac{\prod_{\omega\in\mathcal{N}_1}(k-\omega)}{\prod_{\gamma\in\mathcal{N}_2}(k-\gamma)}\exp\frac{1}{2\pi{\rm i}}\int_{\mathbb{R}}\frac{\log(1+\sigma r_1(s)r_2(s))}{s-k}ds\\
	a_2(k)=\frac{\prod_{\gamma\in\mathcal{N}_2}(k-\gamma)}{\prod_{\omega\in\mathcal{N}_1}(k-\omega)}\exp\frac{1}{2\pi{\rm i}}\int_{\mathbb{R}}\frac{\log(1+\sigma r_1(s)r_2(s))}{s-k}ds
\end{align}
Moreover, in the reflectionless case $r_1(k)=r_2(k)=0$, $a_1, \ a_2$ become
\begin{align}
	a_1(k)=\frac{\prod_{\omega\in\mathcal{N}_1}(k-\omega)}{\prod_{\gamma\in\mathcal{N}_2}(k-\gamma)},\ 	a_2(k)=\frac{\prod_{\gamma\in\mathcal{N}_2}(k-\gamma)}{\prod_{\omega\in\mathcal{N}_1}(k-\omega)}.
\end{align}

With the symmetry condition \eqref{symc}, \eqref{syme},  $M_{sol}(x,t,k)$ satisfies the residue condition \eqref{res1+}, \eqref{res1-} so does $\tilde{M}_{sol}$ where
\begin{align*}
	\tilde{M}_{sol}(x,t;k)=\left\{ \begin{array}{ll}
		\Lambda\overline{M_{sol}(-x,t;-\bar{k})}\Lambda^{-1}\left(\begin{array}{ll}
			\frac{1}{a_1(k)}&0\\
			0&a_1(k)
		\end{array}\right) , &k\in\mathbb{C}^+;\\
		\Lambda\overline{M_{sol}(-x,t;-\bar{k})}\Lambda^{-1}\left(\begin{array}{ll}
			a_2(k)&0\\
			0&\frac{1}{a_2(k)}
		\end{array}\right) , &k\in\mathbb{C}^-.\\
	\end{array}\right.
\end{align*}
That is to say, the $M_{sol}$ satisfies symmetry  \eqref{Msym}.
For convenience, let
\begin{equation}
	M_{sol}=I+\frac{M_{sol}^{(1)}}{k}+\mathcal{O}(\frac{1}{k^2}),\  k\rightarrow \infty.\label{solasy}
\end{equation}
Hence,
\begin{equation}
	M_{sol}^{(1)}=\sum_{\omega_n\in \mathcal{N}_1}
	\left(\begin{array}{cc}
		\alpha_n^{(1)}&0\\
		\alpha_n^{(2)}&0\\
	\end{array}\right)
	+\sum_{\gamma_m\in \mathcal{N}_2}
	\left(\begin{array}{cc}
		0&\beta_m^{(1)}\\
		0&\beta_m^{(2)}\\
	\end{array}\right).\label{msol1}
\end{equation}
For the special case for which all the eigenvalues reside on the imaginary axis, as the  section 9 of \cite{ablo2}, soliton solutions  are considering under $|\mathcal{N}_1|=|\mathcal{N}_2|$. Then
\begin{align}
	a_1(k)=\prod_{j=1}^{|\mathcal{N}_1|}\frac{k-\omega_{j}}{k-\gamma_{j}},\
	a_2(k)=\prod_{j=1}^{|\mathcal{N}_1|}\frac{k-\gamma_{j}}{k-\omega_{j}}.
\end{align}
Letting   $\mathcal{N}_1=\left\{iw|w\in\mathbb{R}^+\right\}$, $\mathcal{N}_2=\left\{-ir|r\in\mathbb{R}^+\right\}$, they derive a breathing one soliton solution
\begin{align}
	q(x,t)=-\frac{2(w+r){\rm e}^{{\rm i}\theta_1}{\rm e}^{-4{\rm i}r^2t}{\rm e}^{-2rx}}{1+{\rm e}^{{\rm i}(\theta_1+\theta_2)}{\rm e}^{4{\rm i}(w^2-r^2)t}{\rm e}^{-2(w+r)x}}.
\end{align}
Where $\theta_1,\theta_2$ are arbitrary real number.

					\subsection{Analysis  on a pure $\bar{\partial}$-Problem}\label{3.4}
					\quad To demonstrate the existence of $M^{RHP}$, we define a new matrix-valued function $W(x,t;k)$
					\begin{equation}
						W(x,t;k)=	M^{(2)}(x,t;k)M^{RHP}(x,t;k)^{-1}.
					\end{equation}
					which   removes   analytic component of $M^{(2)}(z)$    to get  a  pure $\bar{\partial}$-problem.
					According to (\ref{dbar-m2}) we can deduce derivative condition.
					\begin{proposition}\label{prop7}
						$ W(x,t;k)$ satisfies following properties.
						
						$\blacktriangleright$ $ W(x,t;k)$ is continuous   in  $ \mathbb{C}$.
						
						$\blacktriangleright$ $\bar\partial$-equation:\begin{equation}
							\overline{\partial} W=WM^{RHP}\overline{\partial}\mathcal{R} (M^{RHP})^{-1}.
						\end{equation}
						
						$\blacktriangleright$ Asymptotic condition:	\begin{equation}
							W(x,t;k)=I+\mathcal{O}(k^{-1}),\  k\rightarrow \infty,
						\end{equation}
					\end{proposition}
					whose solution can be given by
					\begin{proposition}\label{prop8}
						The Proposition \ref{prop7}  of $W(x,t;k)$ is equivalent to the  integral equation
						\begin{equation}\label{puredbar}
							W(x,t;k)=I-\frac{1}{\pi}\iint_{\mathbb{C}}\frac{W(x,t;s)(M^{RHP}\overline{\partial}\mathcal{R} (M^{RHP})^{-1})(x,t;s)}{s-k}dA(s),
						\end{equation}where $dA(s)$ means integral of general Lebesgue measure on $\mathbb{C}$.
					\end{proposition}
					To solve the (\ref{puredbar}), define $\mathcal{S}$ as the left Cauchy-Green integral  operator with
					\begin{equation}
						\mathcal{S}[f]=-\frac{1}{\pi}\iint_{\mathbb{C}}\frac{f(s)(M^{RHP}\overline{\partial}\mathcal{R} (M^{RHP})^{-1})(x,t;s)}{s-z}dA(s).
					\end{equation}
					Therefore, (\ref{puredbar}) is equivalent to \begin{equation}
						(Id-\mathcal{S})W(x,t;k)=I.
					\end{equation}
					To target on existence of $	W(x,t;k)$, we demonstrate the following proposition.
					\begin{proposition}\label{prop9}
						The norm of the integral operator $\mathcal{S}$ decays to zero as $t\to\infty$:
						\begin{equation}
							||\mathcal{S}||_{\mathcal{B}(L^{\infty})}\lesssim\mathcal{O}(t^{-\frac{1}{4}})+\mathcal{O}(t^{-\frac{1}{2}+\frac{\alpha}{2}}).
						\end{equation}where $\alpha$ is defined in (\ref{alpha}).
					\end{proposition}
					\begin{proof}According to Proposition \ref{dbar-R2},
						\begin{align}
							||\mathcal{S}||_{\mathcal{B}(L^{\infty})}\leq&
							\sum_{j=1}^{4}\frac{1}{\pi}\iint_{\Omega_j}\frac{|M^{RHP}||\overline{\partial}\mathcal{R}||(M^{RHP})^{-1}|}{|s-z|}dA(s)\nonumber\\
							\lesssim&
							\sum_{\Omega_1,\Omega_2,\Omega_3,\Omega_4}\left(\iint_{\Omega_1}\frac{|\overline{\partial}\mathcal{X}|+|r'_1({\rm Re}k)|}{|s-z|}{\rm e}^{-2t|{\rm Re}(i\theta)|}dA(s)\right.\nonumber\\
							&\left.+\iint_{\Omega_1}\frac{|s+\xi|^{-\alpha}}{|s-z|}{\rm e}^{-2t|{\rm Re}(i\theta)|}dA(s)\right). \label{Snorm}
						\end{align}
						Here we have  used the boundedness  of   $M_{RHP}$   in the support of $1-\mathcal{X}$. The first two functions in the last line of above inequalities  are \begin{equation}
							|\overline{\partial}\mathcal{X}|,|r'_1({\rm Re}k)|,|r'_2({\rm Re}k)|,\left|\left( \frac{\sigma r_2}{1+\sigma r_1 r_2}\right) '({\rm Re}k)\right|,\left|\left( \frac{r_1}{1+\sigma r_1 r_2}\right) '({\rm Re}k)\right|.\nonumber
						\end{equation} All of above function are in $L^2(\mathbb{R})$ when the imaginary part of variable is fixed, and $L^2(\mathbb{R})$ norm is independent on the imaginary  part of variable. Therefore, it is valid to use $g(k),||g(u+iv)||_{L_u^2(\mathbb{R})}\leq C_0$ denoting them. In addition, we use $\Omega_1$ as an example to prove
						\begin{align}
							&\iint_{\Omega_1}\frac{|g(s)|}{|s-z|}{\rm e}^{-2t|{\rm Re}(i\theta)|}dA(s) \leq\int_{0}^{\infty}\int_{v}^{\infty}\frac{|g(s)|}{|s-z|}{\rm e}^{-8tuv}dudv\nonumber\\
							&\leq\int_{0}^{\infty}||g||_{L^2_{u}(\mathbb{R})}{\rm e}^{-8tv^2}\left|\left|\frac{1}{(|u|^2+|v-{\rm Im}(z)|^2)^{\frac{1}{2}}}\right|\right|_{L^2_{u}(\mathbb{R})}dv \nonumber\\
							&\lesssim\int_{0}^{\infty}\frac{{\rm e}^{-8tv^2}}{|v-{\rm Im}(z)|^{\frac{1}{2}}}dv
							\leq\int_{\mathbb{R}}\frac{{\rm e}^{-8t(v+{\rm Im}(z))^2}}{|v|^{\frac{1}{2}}}dv=\mathcal{O}(t^{-\frac{1}{4}}).\label{dbarL1}
						\end{align}On the other hand, we demonstrate for $\alpha>0,2>\alpha p>1,p+q=pq$
						\begin{align}
							&\iint_{\Omega_1}\frac{|s+\xi|^{-\alpha}}{|s-z|}{\rm e}^{-2t|{\rm Re}(i\theta)|}dA(s)\nonumber\\
							=&
							\int_{0}^{\infty}\int_{v}^{\infty}\frac{(u^2+v^2)^{-\frac{1}{2}\alpha}{\rm e}^{-8tuv}}{((u-\xi-{\rm Re}(z))^2+(v-{\rm Im}(z))^2)^{\frac{1}{2}}}dudv\nonumber\\
							\leq&\int_{0}^{\infty}{\rm e}^{-8tv^2}\left(\int_{v}^{\infty}(u^2+v^2)^{-\frac{1}{2}\alpha p}du\right)^{\frac{1}{p}}\nonumber\\
							&
							\left(\int_{v}^{\infty}\left((u-\xi-{\rm Re}(z))^2+(v-{\rm Im}(z))^2\right)^{-\frac{1}{2}q}du\right)^{\frac{1}{q}}dv\nonumber\\
							\lesssim&\int_{0}^{\infty}{\rm e}^{-8tv^2}v^{\frac{1}{p}-\alpha}|v-{\rm Im}(z)|^{\frac{1}{q}-1}dv\nonumber	
						\end{align}
						To estimate above integral, we divide the integrating range to two parts:
						\begin{align}	
							&\int_{0}^{\infty}{\rm e}^{-8tv^2}v^{\frac{1}{p}-\alpha}|v-{\rm Im}(z)|^{\frac{1}{q}-1}dv\nonumber\\&=\int_{0}^{{\rm Im}(z)}{\rm e}^{-8tv^2}v^{\frac{1}{p}-\alpha}|v-{\rm Im}(z)|^{\frac{1}{q}-1}dv
							+\int_{{\rm Im}(z)}^{\infty}{\rm e}^{-8tv^2}v^{\frac{1}{p}-\alpha}|v-{\rm Im}(z)|^{\frac{1}{q}-1}dv\nonumber\\
							&={\rm Im}(z)^{1-\alpha}\int_{0}^{1}{\rm e}^{-8t{\rm Im}(z)^2v^2}v^{\frac{1}{p}-\alpha}(1-v)^{\frac{1}{q}-1}dv+\int_{0}^{\infty}{\rm e}^{-8tv^2}v^{\frac{1}{p}-\alpha}(v+{\rm Im}(z))^{\frac{1}{q}-1}dv\nonumber\\
							&\lesssim {\rm Im}(z)^{1-\alpha}\int_{0}^{1}(t^{\frac{1}{2}}{\rm Im}(z)v)^{\alpha-1}v^{\frac{1}{p}-\alpha}(1-v)^{\frac{1}{q}-1}dv+\int_{0}^{\infty}{\rm e}^{-8tv^2}v^{-\alpha}dv\nonumber\\
							&\lesssim\mathcal{O}(t^{-\frac{1}{2}+\frac{\alpha}{2}})+\mathcal{O}(t^{-\frac{1}{2}+\frac{\alpha}{2}})=\mathcal{O}(t^{-\frac{1}{2}+\frac{\alpha}{2}}).
						\end{align}
						According to Proposition \ref{dbar-R2},  substitute $\alpha$ in (\ref{alpha}) into (\ref{Snorm})  together with   (\ref{dbarL1}) yields
						\begin{align}
							||\mathcal{S}||_{\mathcal{B}(L^{\infty})}\lesssim\mathcal{O}(t^{-\frac{1}{4}})+\mathcal{O}(t^{-\frac{1}{2}+\frac{\alpha}{2}})\nonumber	
						\end{align}
					\end{proof}
					Therefore, for sufficiently large $t$, as Proposition \ref{prop7} describing, $W(x,t;k)$ exists. According to Proposition \ref{prop8}, following proposition can be proved.
					
					\begin{proposition}
						$W(x,t;k)$ satisfies large $z$ asymptotic condition
						\begin{align}
							W(x,t;k)&=I+\frac{W^{(1)}(x,t)}{k}+o(k^{-1}),\label{wasy}
						\end{align}
						with
						\begin{align}
							W^{(1)}(x,t)=\frac{1}{\pi}\iint_{\mathbb{C}}W(x,t;s)(M_{RHP}\overline{\partial}\mathcal{R}M_{RHP}^{-1})(x,t;s)dA(s).
						\end{align}
						As $t\to\infty$,
						\begin{align}
							W^{(1)}(x,t)&\lesssim\mathcal{O}(t^{-\frac{3}{4}})+\mathcal{O}(t^{-1+\frac{\alpha}{2}}),
						\end{align}
						where $\alpha$ is defined in (\ref{alpha}).
					\end{proposition}
					\begin{proof}
						As Proposition \ref{prop8} shows, $W(x,t;k)$ has integral form \begin{equation}
							W(x,t;k)=I-\frac{1}{\pi}\iint_{\mathbb{C}}\frac{W(x,t;s)(M^{RHP}\overline{\partial}\mathcal{R} (M^{RHP})^{-1})(x,t;s)}{s-k}dA(s).\nonumber
						\end{equation}
						To compute the large $k$ asymptotic property, we write $W(x,t;k)$ into
						\begin{align*}
							W(x,t;k)=&I+k^{-1}\left[\frac{1}{\pi}\iint_{\mathbb{C}}W(x,t;s)(M^{RHP}\overline{\partial}\mathcal{R} (M^{RHP})^{-1})(x,t;s)dA(s)\right]\\
							&-\frac{1}{\pi}\iint_{\mathbb{C}}\frac{sW(x,t;s)(M^{RHP}\overline{\partial}\mathcal{R} (M^{RHP})^{-1})(x,t;s)}{k(s-k)}dA(s).
						\end{align*}
						Similar to \cite{dbar1}, (\ref{wasy}) is proved. According to Proposition \ref{dbar-R2}, we have
						\begin{align*}
							W^{(1)}(x,t)=&\frac{1}{\pi}\iint_{\mathbb{C}}W(x,t;s)(M^{RHP}\overline{\partial}\mathcal{R} (M^{RHP})^{-1})(x,t;s)dA(s)\\
							\lesssim&||W||_{L^{\infty}}\sum_{j=1}^{4}\iint_{\Omega_j}|M^{RHP}||\overline{\partial}\mathcal{R} (M^{RHP})^{-1}|dA(s)\nonumber\\
							\lesssim&
							\sum_{\Omega_1,\Omega_2,\Omega_3,\Omega_4}\left[\iint_{\Omega_1}\left( |\overline{\partial}\mathcal{X}|+|\overline{\partial}\mathcal{Y}|+|r'_1({\rm Re}k)|\right) {\rm e}^{-2t|{\rm Re}({\rm i}\theta)|}dA(s)\right.\nonumber\\
							&\left.+\iint_{\Omega_1}\left( |s+\xi|^{-\frac{1}{2}-{\rm Im}\nu(-\xi)}+|s+\xi|^{-\frac{1}{2}}\right) {\rm e}^{-2t|{\rm Re}({\rm i}\theta)|}dA(s)\right].
						\end{align*}
						Similar with the proof of Proposition \ref{prop9}, for $g(k),||g(u+iv)||_{L_u^2(\mathbb{R})}\leq C_0$, there is
						\begin{align}
							\iint_{\Omega_1}|g(s)| {\rm e}^{-2t|{\rm Re}({\rm i}\theta)|}dA(s)\leq&\int_{0}^{\infty}\int_{v}^{\infty}|g(s)|{\rm e}^{-8tuv}dudv\nonumber\\
							\leq&\int_{0}^{\infty}||g(u+{\rm i}v)||_{L_u^2(\mathbb{R})}\int_{v}^{\infty}|g(s)|{\rm e}^{-16tuv}dudv\nonumber\\
							\lesssim&t^{-\frac{1}{2}}\int_{0}^{\infty}\frac{{\rm e}^{-8tv^2}}{\sqrt{v}}dv
							\lesssim \mathcal{O}(t^{-\frac{3}{4}}).\nonumber
						\end{align}
						Meanwhile, for $ \alpha>0,\ 2>\alpha p>1,\ p+q=pq$
						\begin{align*}
							&\iint_{\Omega_1}|s+\xi|^{-\alpha}{\rm e}^{-2t|{\rm Re}({\rm i}\theta)|}dA(s)\nonumber
							\leq
							\int_{0}^{\infty}\int_{v}^{\infty}(u^2+v^2)^{-\frac{1}{2}\alpha}{\rm e}^{-8tuv}dudv\nonumber\\
							\leq&\int_{0}^{\infty}\left(\int_{v}^{\infty}{\rm e}^{-8qtuv}du\right)^{\frac{1}{q}}\left(\int_{v}^{\infty}(u^2+v^2)^{-\frac{1}{2}\alpha p}du\right)^{\frac{1}{p}}dv\nonumber\\
							\lesssim&\int_{0}^{\infty}(tv)^{-\frac{1}{q}}{\rm e}^{-8tv^2}v^{\frac{1}{p}-\alpha}dv\\
							=&t^{-\frac{1}{q}}\int_{0}^{\infty}v^{\frac{1}{p}-\frac{1}{q}-\alpha}{\rm e}^{-8tv^2}dv = \mathcal{O}(t^{-1+\frac{\alpha}{2}}).
						\end{align*}
						By substituting $\alpha$ in (\ref{alpha}) into the above result, we deduce
						\begin{align}
							W^{(1)}(x,t)\lesssim\mathcal{O}(t^{-\frac{3}{4}})+\mathcal{O}(t^{-1+\frac{\alpha}{2}}).
						\end{align}
					\end{proof}

					\section{Long time asymptotic behavior}\label{sec4}
					
					By summarizing the transformation above, we have
					\begin{equation}
						M(x,t;k)=W(x,t;k)E(x,t;k)M_{sol}(x,t;k)\mathcal{R}(x,t;k)\delta^{\sigma_3}(k).
					\end{equation}
					From (\ref{solasy}), (\ref{Easy}), (\ref{wasy}), and  $\mathcal{R}=I$ on the  imaginary axis, large $k$ along the imaginary axis expansion of $M(x,t;k)$ can be deduced as
					\begin{align*}
						M(x,t;k)=&\left(I+\frac{W^{(1)}(x,t)}{k}+o(k^{-1})\right) \left(I+\frac{E^{(1)}(x,t)}{k}+\mathcal{O}(k^{-2})\right)\\
						&\left(I+\frac{M_{sol}^{(1)}}{k}+\mathcal{O}(k^{-2})\right)
						\left(I+\frac{\delta^{(1)}\sigma_3}{k}+\mathcal{O}(k^{-2})\right)\\
						=&I+k^{-1}\left(W^{(1)}(x,t)+E^{(1)}(x,t)+M_{sol}^{(1)}+\delta^{(1)}\sigma_3\right)+o(k^{-1})
					\end{align*}
					Denote $M^{(1)}(x,t)$ as the coefficient of $1/k$ term in above expansion of $M(x,t;k)$, then
					\begin{align*}
						M^{(1)}(x,t)=&W^{(1)}(x,t)+E^{(1)}(x,t)+M_{sol}^{(1)}+\delta^{(1)}\sigma_3\\
						=&M_{sol}^{(1)}+\delta^{(1)}\sigma_3+\frac{1}{\sqrt{8t}}\left(\begin{array}{cc}
							0&-{\rm i}t^{{\rm Im}\nu(-\xi)}\tilde{\beta}_{12}\\
							{\rm i}t^{-{\rm Im}\nu(-\xi)}\tilde{\beta}_{21}&0\\
						\end{array}\right)\\
						+&\mathcal{O}(t^{-\frac{1}{2}+|{\rm Im}\nu(-\xi)|})\begin{array}{cc}
							\Big(\mathcal{O}(t^{-\frac{1}{2}-{\rm Im}\nu(-\xi)}),\mathcal{O}(t^{-\frac{1}{2}+{\rm Im}\nu(-\xi)})\Big)\end{array}\\
						+&\mathcal{O}(t^{-\frac{3}{4}})+\mathcal{O}(t^{-1+\frac{\alpha}{2}}).
					\end{align*}
					Using the (\ref{qfml}), we can reconstruct a  solution of (\ref{NNLS}) denote as $q_{sol}$ where the modified  scattering data
					\begin{align}
							\Big\lbrace r_1\equiv0,r_2\equiv0,\mathcal{N}_1,\mathcal{N}_2,\mathcal{P}_1',\mathcal{P}_2'\Big\rbrace,
						\end{align} where \begin{align*}
							\mathcal{P}_1'=&\left\{\delta^{-2}(\omega_{n})c_n,\delta^{-2}(\eta_{m})e_m|n=\pm1,..\pm N_1,m=1,..M_1\right\},\\
							\mathcal{P}_2'=&\left\{\delta^2(\gamma_{n})d_n,\delta^{2}(\tau_{m})f_m|n=\pm1,..\pm N_2,m=1,..M_2\right\}.
					\end{align*}
					On this occasion, the solution of RHP0 is $M_{sol}$. Therefore, according to \eqref{msol1}
					\begin{align}\label{qsol}
						q_{sol}=2{\rm i}\lim_{k\rightarrow\infty}\left(k[M_{sol}(x,t;k)]_{12}\right)=2{\rm i}\sum_{\gamma\in \mathcal{N}_2}\beta_\gamma^{(1)},
					\end{align} where $\beta_\gamma$ is solution of \eqref{res1},\eqref{res2}.
					Then we achieve main result of this paper.
					\begin{theorem}\label{thm1} Let $q(x,t)$ be the solution for  the initial-value problem (\ref{NNLS}) with generic data  $q_0 \in H^{1,1}$ such that  assumption \ref{asmp1} and assumption \ref{asym2}
						are satisfied.
						Let $\xi=\frac{x}{4t}\in\mathbb{R}$,
						there exists a large constant $T_1=T_1(\xi)$, for all $t>T_1$, long-time asymptotic  behavior of the  solution $q(x,t)$ is given as follows
						\begin{align}\label{main}
							q(x,t)= q_{sol}(x,t)+\frac{t^{{\rm Im}\nu(-\xi)}\tilde{\beta}_{12}}{\sqrt{2t}}+\left\lbrace\begin{array}{llll}
								\mathcal{O}(t^{-1+2{\rm Im}\nu(-\xi)}),&{\rm Im}\nu(-\xi)\in(\frac{1}{6},\frac{1}{2})\\[12pt]
								\mathcal{O}(t^{-\frac{3}{4}+\frac{{\rm Im}\nu(-\xi)}{2}}),&{\rm Im}\nu(-\xi)\in(0,\frac{1}{6}]\\[12pt]
								\mathcal{O}(t^{-\frac{3}{4}}),&{\rm Im}\nu(-\xi)\in(-\frac{1}{4},0]
							\end{array}\right.,
						\end{align}
						
						and \begin{align}
							q(x,t)= q_{sol}(x,t)+
							\mathcal{O}(t^{-\frac{3}{4}}),\hspace{0.5cm}{\rm Im}\nu(-\xi)\in(-\frac{1}{2},-\frac{1}{4}].
						\end{align}
						Above $q_{sol}$, $\nu$ and $\tilde{\beta}_{12}$ are defined by (\ref{qsol}), (\ref{nu}) and  (\ref{b12}), respectively.
					\end{theorem}


\Acknowledgements{This work is supported by  the National Natural Science
	Foundation of China (Grant No. 12271104, 51879045, 12247182). The authors
	would like to have their sincerest gratitude to referees for patient
	guidance and valuable suggestions.}





\end{document}